\documentclass[11pt,a4paper]{article}

\usepackage[T1]{fontenc}
\usepackage{lmodern}
\usepackage[margin=0.9in]{geometry}
\usepackage{amsmath,amssymb,amsthm}
\usepackage{microtype}
\usepackage[hidelinks]{hyperref}

\newtheorem{theorem}{Theorem}

\newcommand{\C}{\mathbb{C}}
\newcommand{\T}{\mathbb{T}}

\title{A proof of the constant-width numerical range conjecture}
\author{
  Muyan Jiang\\[0.3em]
  \small University of California, Berkeley\\
  \small \href{mailto:muyan_jiang@berkeley.edu}{muyan\_jiang@berkeley.edu}
  \and
  Minwu Kim\\[0.3em]
  \small Seoul National University\\
  \small \href{mailto:minwukim@snu.ac.kr}{minwukim@snu.ac.kr}
}
\date{}

\begin{document}

\maketitle

\begin{abstract}
We prove that the numerical range of a complex matrix has constant width \(d\) if and only if it is a closed disk of radius \(d/2\). The width in each direction is the gap between the extreme eigenvalues of a Hermitian matrix. When that gap is constant, the midpoint of the projected interval is a rational function of the direction parameter. An operator norm bound and reversal of direction then show that the midpoints are projections of a single point, the disk's center.
\end{abstract}

The numerical range of a complex matrix \(A\in\C^{n\times n}\) is
\[
  W(A)=\{x^*Ax:x\in\C^n,\ x^*x=1\},
\]
where \({}^*\) denotes conjugate transpose. The set \(W(A)\) is compact
and, by the Toeplitz--Hausdorff theorem \cite{Hausdorff1919}, convex.
Its width in a given direction is the length of its orthogonal
projection onto the line through the origin in that direction. The Reuleaux
triangle is a familiar example of a noncircular set of constant width.
Chien, Nakazato, and Meng \cite{CNM2019} conjectured that a numerical
range of constant width must be a disk. They studied constant-width
examples among nilpotent Toeplitz matrices. Chien, Li, and Nakazato
\cite[Sec.~1]{CLN2021} later proposed the analogous statement for
higher-rank numerical ranges.

\begin{theorem}
For any \(A\in\C^{n\times n}\), \(n\geq1\), and \(d\geq0\), the
numerical range \(W(A)\) has constant width \(d\) if and only if it is a
closed disk of radius \(d/2\). A disk of radius zero is a singleton.
\end{theorem}

\begin{proof}
A disk of radius \(d/2\) has constant width \(d\). Now suppose \(W(A)\)
has constant width \(d\). If \(d=0\), the projections onto the real and
imaginary axes are both points, so \(W(A)\) is a singleton. For the rest
of the proof let \(r=d/2>0\).

Write \(\T\) for the unit circle and define
\[
  H(\zeta)=\frac{\zeta^{-1}A+\zeta A^*}{2},
  \qquad \zeta\in\C\setminus\{0\}.
\]
For \(\zeta\in\T\), the projection in direction \(\zeta\) has coordinates
\(\operatorname{Re}(\bar\zeta w)\) for \(w\in W(A)\). Since
\(\zeta^{-1}=\bar\zeta\), the matrix \(H(\zeta)\) is Hermitian and
\(x^*H(\zeta)x=\operatorname{Re}(\bar\zeta x^*Ax)\). By Rayleigh--Ritz
\cite{HornJohnson2013}, the endpoints of the projected interval are the
smallest and largest eigenvalues of \(H(\zeta)\). Let \(m(\zeta)\) be
the midpoint, so the endpoints are \(m(\zeta)\pm r\). We seek a fixed
\(c\in\C\) with \(m(\zeta)=\operatorname{Re}(\bar\zeta c)\) for all
\(\zeta\in\T\), as for a disk centered at \(c\).

In a fixed direction \(\zeta_0\in\T\), two eigenvalues separated by the
full width \(2r\) must be the endpoints. Thus \(m(\zeta_0)\) is the only
\(t\in\C\) for which both \(t-r\) and \(t+r\) are eigenvalues of
\(H(\zeta_0)\). With \(I_n\) denoting the \(n\times n\) identity matrix,
form the shifted characteristic polynomials
\[
  P_\pm(t,\zeta)=\det\bigl((t\pm r)I_n-H(\zeta)\bigr).
\]

Regard \(P_-\) and \(P_+\) as polynomials in \(t\) over the field
\(\C(\zeta)\) of rational functions in \(\zeta\). If they were coprime,
B\'ezout's identity would give \(UP_-+VP_+=1\) for some
\(U,V\in\C(\zeta)[t]\). Choose \(\zeta_0\in\T\) away from the finitely
many poles of the coefficients of \(U\) and \(V\). At \(\zeta=\zeta_0\)
and \(t=m(\zeta_0)\) both shifted polynomials vanish, giving \(0=1\), a contradiction.

Let \(g\) be their monic greatest common divisor, of degree \(k\geq1\)
in \(t\). At each \(\zeta_0\in\T\) where the coefficients of \(g\) and
both quotients \(P_\pm/g\) are finite, \(g(t,\zeta_0)\) is still monic
of degree \(k\) and divides both evaluated polynomials. Their only
common root is \(m(\zeta_0)\), so \(g(t,\zeta_0)=(t-m(\zeta_0))^k\).
Define \(f\in\C(\zeta)\) by writing the coefficient of \(t^{k-1}\) in
\(g\) as \(-kf(\zeta)\). Then \(f(\zeta_0)=m(\zeta_0)\) away from the
finitely many excluded points. Both rational functions
\(P_\pm(f(\zeta),\zeta)\) therefore have infinitely many zeros, so
\begin{equation}\label{eq:characteristic-identities}
  P_-(f(\zeta),\zeta)=P_+(f(\zeta),\zeta)=0
  \qquad\text{in }\C(\zeta).
\end{equation}

For any nonzero complex \(\zeta\) where \(f\) is finite,
\eqref{eq:characteristic-identities} gives \(f(\zeta)+r\) as an eigenvalue
of \(H(\zeta)\). In the Euclidean operator norm \(\|\cdot\|\), we have
\(|f(\zeta)+r|\leq\|H(\zeta)\|\) and \(\|A^*\|=\|A\|\). Hence
\[
  |\zeta f(\zeta)|
  \leq |\zeta|\bigl(\|H(\zeta)\|+r\bigr)
  \leq \frac{\|A\|}{2}(1+|\zeta|^2)+r|\zeta|.
\]
This bound rules out a pole of \(\zeta f(\zeta)\) at any finite point,
including zero. Thus \(\zeta f(\zeta)\) is a polynomial. It grows at most
quadratically at infinity, so its degree is at most two and
\[
  f(\zeta)=a_{-1}\zeta^{-1}+a_0+a_1\zeta,
  \qquad a_{-1},a_0,a_1\in\C.
\]
In particular, \(f\) is finite on the whole unit circle. Evaluating
\eqref{eq:characteristic-identities} there gives \(f=m\) in every
direction by uniqueness of the common root.

Reversing the direction negates the midpoint, so
\(f(-\zeta)=-f(\zeta)\) on \(\T\) and hence \(a_0=0\). The function
\(f=m\) is real on \(\T\), where \(\bar\zeta=\zeta^{-1}\). Multiplying
\(f(\zeta)=\overline{f(\zeta)}\) by \(\zeta\) gives an equality of
polynomials on \(\T\), and comparing coefficients yields
\(a_1=\overline{a_{-1}}\). Thus, with \(c=2a_{-1}\),
\[
  m(\zeta)=\operatorname{Re}(\bar\zeta c),
  \qquad \zeta\in\T.
\]
The upper endpoint \(\operatorname{Re}(\bar\zeta c)+r\) is the support
function of both \(W(A)\) and the closed disk centered at \(c\) with
radius \(r\). A compact convex set is determined by its support function
\cite[Sec.~1.7]{Schneider2014}, so the two sets coincide.
\end{proof}

\section*{Acknowledgments}

The authors thank Ilya Spitkovsky for suggesting the use of
characteristic roots separated by the width and for his comments
on earlier drafts.

The authors used OpenAI's ChatGPT to explore proof strategies
and to develop and refine the proof and its exposition.
The authors take responsibility for the mathematical arguments
and the final manuscript.


\begin{thebibliography}{9}
\raggedright

\bibitem{Hausdorff1919}
F. Hausdorff,
\emph{Der Wertvorrat einer Bilinearform},
Math. Z. \textbf{3} (1919), 314--316.
\href{https://doi.org/10.1007/BF01292610}
{doi:10.1007/BF01292610}.

\bibitem{CNM2019}
M.-T. Chien, H. Nakazato, and J. Meng,
\emph{The diameter and width of numerical ranges},
Linear Algebra Appl. \textbf{582} (2019), 76--98.
\href{https://doi.org/10.1016/j.laa.2019.07.036}
{doi:10.1016/j.laa.2019.07.036}.

\bibitem{CLN2021}
M.-T. Chien, C.-K. Li, and H. Nakazato,
\emph{The diameter and width of higher rank numerical ranges},
Linear Multilinear Algebra \textbf{69} (2021), no.~5, 871--887.
\href{https://doi.org/10.1080/03081087.2019.1710105}
{doi:10.1080/03081087.2019.1710105}.

\bibitem{HornJohnson2013}
R. A. Horn and C. R. Johnson,
\emph{Matrix Analysis}, 2nd ed.,
Cambridge University Press, Cambridge, 2013.
\href{https://doi.org/10.1017/CBO9781139020411}
{doi:10.1017/CBO9781139020411}.

\bibitem{Schneider2014}
R. Schneider,
\emph{Convex Bodies: The Brunn--Minkowski Theory}, 2nd expanded ed.,
Encyclopedia of Mathematics and its Applications, vol.~151,
Cambridge University Press, Cambridge, 2014.
\href{https://doi.org/10.1017/CBO9781139003858}
{doi:10.1017/CBO9781139003858}.

\end{thebibliography}
\end{document}